\documentclass[a4paper,12pt,reqno]{amsart}
\usepackage{verbatim}
\usepackage{enumitem}
\usepackage{amssymb}
\usepackage[all]{xy}
\usepackage{mathtools}
\usepackage{color}
\usepackage{tabu}
\usepackage{longtable}
\usepackage{float}
\usepackage{amsthm}
\usepackage{booktabs}
\definecolor{cobalt}{RGB}{61,89,171}
\usepackage[colorlinks,citecolor=cobalt,linkcolor=cobalt,urlcolor=cobalt,pdfpagemode=UseNone,backref = page]{hyperref}

\newcommand{\ox}{\otimes}
\newcommand{\RR}{\mathbb{R}}
\newcommand{\la}{\langle}
\newcommand{\ra}{\rangle}

\newcommand{\G}{\Gamma}
\newcommand{\vol}{\mathrm{vol}}

\newcommand{\R}{\mathbb{R}}
\newcommand{\lie}[1]{\mathfrak{#1}}     

\newcommand{\N}{\mathbb{N}}

\newcommand{\SU}{\mathrm{SU}}

\newcommand{\Gtwo}{\mathrm{G}_2}

\newcommand{\fg}{\mathfrak{g}}
\newcommand{\fh}{\mathfrak{h}}
\newcommand{\fn}{\mathfrak{n}}

\newtheorem{theorem}{Theorem}[section]
\newtheorem*{theorem*}{Theorem}

\newtheorem{proposition}[theorem]{Proposition}

\newtheorem{example}[theorem]{Example}

\newtheorem{question}[theorem]{Question}

\theoremstyle{remark}
\newtheorem{remark}[theorem]{Remark}

\begin{document}

\title{Purely coclosed $\mathbf{G}_2$-structures on nilmanifolds - II}

\author[G. Bazzoni]{Giovanni Bazzoni}
\address{Dipartimento di Scienza ed Alta Tecnologia, Universit\`a degli Studi dell'Insubria, Via Valleggio 11, 22100, Como, Italy}
\email{giovanni.bazzoni@uninsubria.it}

\author[G. Petracci]{Giorgia Petracci}
\address{Dipartimento di Scienza ed Alta Tecnologia, Universit\`a degli Studi dell'Insubria, Via Valleggio 11, 22100, Como, Italy}
\email{gpetracci@studenti.uninsubria.it}




\begin{abstract}
This paper completes the classification of seven-dimensional nilpotent Lie groups endowed with a left-invariant purely coclosed $\text{G}_2$-structure, initiated in \cite{purely1}. In \cite{purely1}, the authors provided the classification of decomposable seven-dimensional nilpotent Lie groups and of the indecomposable ones up to step $4$ of nilpotency. Here, we address the case of indecomposable $5$- and $6$-step nilpotent Lie groups.
\end{abstract}

\maketitle

\section{Introduction}

$\Gtwo$ is one of the exceptional Lie groups appearing in the Killing-Cartan classification of semisimple Lie groups. Aside from its relevance in the theory of Lie groups, $\Gtwo$ plays a crucial role in Riemannian geometry. In this context, its importance stems from a celebrated result of Berger, stating that $\Gtwo$ is the only (reduced) holonomy group of an irreducible, non locally symmetric Riemannian 7-manifold \cite{Berger}. Such Riemannian manifolds are in particular Ricci-flat; in fact, the only known examples of odd-dimensional, compact Riemannian Ricci-flat manifolds are provided by 7-manifolds with holonomy (contained in) $\Gtwo$. Notice that Berger's result says that $\Gtwo$ is a possible holonomy group. In fact, the construction of manifolds with holonomy $\Gtwo$ has proved to be a very difficult task, see \cite{Bryant-Salamon,Joyce,Joyce-Karigiannis}. Such manifolds also play a relevant role in $M$-theory, see for instance \cite{Atiyah-Witten}.

Leaving the holonomy context, a $\Gtwo$-structure on a $7$-dimensional manifold is a reduction of its frame bundle from $\mathrm{GL}(7,\RR)$ to $\Gtwo \subset \mathrm{SO}(7)$; as shown by Gray \cite{Gray}, the existence of such a reduction is equivalent to the manifold being orientable and spin. Thus a $\Gtwo$-structure is equivalent to the existence of a positive $3$-form $\varphi$ (see \cite{Bryant} for details). This form defines a Riemannian metric $g_{\varphi}$ and an orientation $\vol_\varphi$ on $M$, hence a Hodge operator $*_\varphi$. When $\varphi$ is parallel with respect to the Levi-Civita connection of $g_\varphi$, the identity component of the holonomy of $g_\varphi$ is contained in $\Gtwo$, and we are in the situation of the last paragraph. By a result of Fern\'andez and Gray, $\varphi$ is parallel if and only if $\varphi$ is closed and coclosed, see \cite{FernandezGray}. In this case, $g_\varphi$ is Ricci-flat. A $\Gtwo$-structure is called {\it closed} if $d\varphi=0$, and {\it coclosed} if $d*_\varphi\varphi=0$. By a recent result of Crowley and Nordstr\"om, coclosed $\Gtwo$-structures exist on any closed, oriented spin manifold, since they satisfy an $h$-principle (see \cite[Theorem 1.8]{CN}).

As it is the case for general $\mathrm{G}$-structures, the non-integrability of a $\Gtwo$-structure is governed by its intrinsic torsion $\tau$, see \cite{Salamon}. In this particular case,
$\tau$ has four components $\tau_i\in\Omega^i(M)$, $i=0,1,2,3$, determined by the equations
\begin{equation}\label{Bryant_torsion}
\left\{
\begin{array}{lcl}
d\varphi & = & \tau_0*_\varphi\varphi + 3~\tau_1\wedge \varphi + *_\varphi\tau_3 \\
d*_\varphi\varphi & = & 4\tau_1 \wedge *_\varphi\varphi + \tau_2\wedge \varphi 
\end{array}
\right.\,;
\end{equation}
see \cite[Proposition 1]{Bryant}. Closed $\Gtwo$-structures are thus those for which $\tau_0=\tau_1=\tau_3=0$, while coclosed $\Gtwo$-structures are characterized by $\tau_1=\tau_2=0$. A $\Gtwo$-structure is of {\em pure type} if all the torsion components vanish, but one. Thus closed $\Gtwo$-structures are of pure type, while coclosed $\Gtwo$-structures are not. $\Gtwo$-structures with  $\tau_i=0$, $i=0,2,3$, are called {\em locally conformally parallel} (see \cite{IPP}). $\Gtwo$-structures with $\tau_i=0$, $i=1,2,3$ are known as {\em nearly parallel}, see \cite{FKMS}. In this case $d\varphi=\tau_0*_\varphi\varphi$, $\tau_0$ is a constant, and $g_\varphi$ is Einstein with positive scalar curvature. Nearly parallel $\Gtwo$-structures are in particular coclosed.

This paper deals with the last pure class of $\Gtwo$-structures, namely {\em purely coclosed $\Gtwo$-structures}; these are 
given by the conditions $\tau_0=\tau_1=\tau_2=0$, that is, $d\varphi=*_\varphi\tau_3$ and $d*_\varphi\varphi=0$. Equivalently (see \cite{BMR}), they can be characterized by
\[
d*_\varphi\varphi=0 \quad \mathrm{and} \quad \varphi\wedge d\varphi=0\,.
\]
In particular, a purely coclosed $\Gtwo$-structure is coclosed. The second equation is an equality of $7$-forms, so it imposes a single extra condition. It is not clear whether there is an $h$-principle for purely coclosed $\Gtwo$-structures.

We study left-invariant $\Gtwo$-structures on nilmanifolds, compact quotients of connected, simply connected, nilpotent Lie groups by a lattice. Hence, we can restrict the attention to $7$-dimensional real nilpotent Lie algebras, which have been classified by Gong in \cite{Gong}. As for pure classes, nilpotent Lie groups endowed with a left-invariant closed $\Gtwo$-structure have been classified by Conti and Fern\'andez in \cite{CF}. It can be shown that nilmanifolds can not carry locally conformally parallel or left-invariant nearly parallel $\Gtwo$-structures (see the discussion in \cite{purely1}). As for left-invariant (purely) coclosed $\Gtwo$-structures, we report some known facts from the literature:
\begin{itemize}
    \item Bagaglini, Fern\'andez and Fino determined in \cite{BFF} which nilpotent Lie groups admit left-invariant coclosed $\Gtwo$-structures in two cases: when the Lie algebra is decomposable, and when it is $2$-step nilpotent.
    \item In \cite{Freibert}, Freibert obtained all nilpotent almost-abelian Lie algebras with a coclosed $\Gtwo$-structure.
    \item del Barco, Moroianu and Raffero classified $2$-step nilpotent Lie groups admitting left-invariant purely coclosed invariant $\Gtwo$-structures. Their approach uses the theory of 2-step nilpotent Lie algebras and does not rely on the classification of $7$-dimensional nilpotent Lie algebras (see \cite{BMR}).
    \item In \cite{purely1}, the first author, together with Garvín and Mu\~noz, determined which decomposable and indecomposable Lie algebras of nilpotency step $\leq 4$ admit purely coclosed $\Gtwo$-structures.
\end{itemize}

We refer the reader to Section \ref{sec:preliminaries} for the relevant details.

In this paper we follow the same approach as in \cite{purely1} and determine which 7-dimensional indecomposable nilpotent Lie groups with nilpotency step $\geq 5$ admit a left-invariant purely coclosed $\Gtwo$-structure. In the positive cases, we provide an explicit example of a left-invariant purely coclosed $\Gtwo$-structure. In the negative cases, we consider obstructions to the existence of coclosed $\Gtwo$-structures. The results are summarized in Theorems \ref{5step-indecomposable-pure} and \ref{6step-indecomposable-pure}. In particular, putting together our results with those of \cite{purely1} and \cite{BMR}, we have (see Theorem \ref{thm:nice}):
\begin{theorem*}\label{thm:main}
Let $\fg$ be a 7-dimensional nilpotent Lie algebra.
\begin{itemize}
    \item If $\fg$ is decomposable, then it admits a coclosed $\Gtwo$-structure if and only if it admits a purely coclosed one, unless $\fg=\fh_3\oplus\R^4$.
    \item If $\fg$ is indecomposable, then it admits a coclosed $\Gtwo$-structure if and only if it admits a purely coclosed one.
\end{itemize}
\end{theorem*}
Here, $\fh_3$ denotes the Heisenberg Lie algebra.


The computations are implemented using \texttt{SageMath} \cite{Sage}. For each of the Lie algebras there is a single worksheet. When the Lie algebra admits a left-invariant purely coclosed $\Gtwo$-structure, we provide the relevant forms, and the list of commands that verify that they define a purely coclosed $\Gtwo$-structure. In the negative case,
we provide a commented worksheet with all the steps that show the obstructions. The SageMath worksheets can be found at \cite{worksheets}.

\noindent {\bf Acknowledgements.} 
We thank Luc\'ia Mart\'in-Merch\'an and Vicente Mu\~noz for useful discussions. We also thank the anonymous referees and the editor for their valuable comments, which greatly improved the exposition. The first author is partially supported by the PRIN 2022 project ``Interactions between Geometric Structures and Function Theories''(code 2022MWPMAB) and by the GNSAGA of INdAM.

\section{Preliminaries}\label{sec:preliminaries}

This work uses the same techniques introduced in \cite {purely1} in order to finish the classification of 7-dimensional nilpotent Lie algebras with purely coclosed $\Gtwo$-structures. Nevertheless, we take the opportunity to include a short section with the necessary preliminaries, in order to make this paper as self-contained as possible.

\subsection{$\mathbf{G}_2$- and $\mathbf{SU}(3)$-structures}

In this paper we describe $\Gtwo$-structures using a positive 3-form. More precisely, a $\Gtwo$-structure on a 7-dimensional manifold $M$ is a 3-form $\varphi\in\Omega^3(M)$ for which at each point there exists a local coframe $\{e^i\}$ such that 
\begin{equation}\label{G2-form}
\varphi = e^{127}+e^{347}+e^{567}+e^{135}-e^{146}-e^{236}-e^{245}\,.    
\end{equation}
Such 3-forms are called {\em positive} and the space of positive 3-forms is denoted $\Omega^3_+(M)$. Then $\varphi$ defines a Riemannian metric $g_\varphi$ and a volume form $\mathrm{vol}_\varphi$ by the formula
\[
g_{\varphi}(X,Y)\, \mathrm{vol}_{\varphi} = \frac{1}{6}\bigl(\imath_{X}\varphi\bigr)\wedge\bigl(\imath_{Y}\varphi\bigr)\wedge \varphi\,,
\]
for vector field $X,Y$ on $M$. General references on $\Gtwo$-geometry are \cite{Bryant,Karigiannis}. We focus on $\Gtwo$-structures of {\em purely coclosed} type: with respect to \eqref{Bryant_torsion}, they are defined by the vanishing of all torsion forms but $\tau_3$. Equivalently, they can be defined as coclosed $\Gtwo$-structures, i.e. $d\ast_\varphi\varphi=0$, such that $\varphi\wedge d\varphi=0$.

We are interested in particular instances of $\Gtwo$-structures on nilmanifolds, which are compact quotients of connected, simply connected, nilpotent Lie groups. Let $G$ be a $7$-dimensional connected and simply connected Lie group with Lie algebra $\lie{g}$. A $\Gtwo$-structure on $G$ is \emph{left-invariant} if the defining $3$-form is left-invariant. Thus, a left-invariant $\Gtwo$-structure on $G$ is defined by a positive $3$-form $\varphi\in \Lambda^3{\lie{g}}^*$ which can be written, in some orthonormal coframe $\{e^1,\dotsc, e^7\}$ of ${\mathfrak{g}}^*$, as \eqref{G2-form}. A $\Gtwo$-structure on $\lie{g}$ is {\em coclosed} if $\varphi$ is coclosed, that is, if
\[
d*_\varphi\varphi=0\,,
\]
where $d$ denotes the Chevalley-Eilenberg differential on ${\lie{g}}^*$. A coclosed $\Gtwo$-structure on $\lie{g}$ is said to be {\em purely coclosed} if, in addition, 
\[
\varphi\wedge d\varphi=0\,.
\]

$\Gtwo$-structures on 7-dimensional manifolds are strictly related to $\mathrm{SU}(3)$-structures on $6$-manifolds. Indeed, both structures can be described uniformly in terms of a certain spinor, see \cite{ACFH}. On a 6-dimensional manifold $N$ consider a pair $(\psi_-,\omega)\in\Omega^3(N)\times\Omega^2(N)$, where $\psi_-$ is stable in the sense of Hitchin and $\omega$ is almost symplectic. There are two orbits of stable forms in $\Omega^3(N)$: the {\em positive} and the {\em negative} stable forms (see \cite[Proposition 2]{Hitchin}). If $\psi_-$ is negative, one can use it to define an almost complex structure $J$ on $N$. We call $(\psi_-,\omega)$ an $\mathrm{SU}(3)$-structure on $N$ if, in addition, $\psi_-\wedge\omega=0$ and $g(\cdot,\cdot)=\omega(\cdot,J\cdot)$ defines a Riemannian metric on $N$. In this case, defining $\psi_+=-J^*\psi_-$, one obtains a complex volume form $\Psi=\psi_++i\psi_-$ on $N$.

Starting with an $\mathrm{SU}(3)$-structure $(\psi_-,\omega)$ on a 6-dimensional manifold $N$, we can construct a $\Gtwo$-structure on $M=N\times\RR$ by setting $\varphi=\omega\wedge dt+\psi_+$, where $t$ is a coordinate on $\RR$.

As it was the case for $\Gtwo$-structures, one can make sense of left-invariant $\mathrm{SU}(3)$-structures on Lie groups.

\subsection{Nilmanifolds} Given a Lie algebra $\fg$, consider the series of ideals defined by
\[
\fg_0=\{0\}\,, \qquad \fg_k=\{X\in\fg \mid [X,\fg]\subset\fg_{k-1}\} \quad \mathrm{for} \ k\geq 1\,.
\]
Notice that $\fg_1=\mathfrak{z}(\fg)$, the center of $\fg$. $\fg$ is {\em nilpotent} if there exists $n\in\N$ such that $\fg_n=\fg$. In particular, $\mathfrak{z}(\fg)\neq 0$ if $\fg$ is nilpotent. Moreover, $\fg$ is {\em n-step nilpotent} if $\fg_n=\fg$ but $\fg_{n-1}\neq \fg$; $n$ is then the {\em nilpotency step} of $\fg$. A connected Lie group $G$ is {\em nilpotent} or {\em n-step nilpotent} if its Lie algebra $\fg$ is.

A nilmanifold $M=\G\backslash G$ is the compact quotient of a connected, simply connected, nilpotent Lie group $G$ by a lattice $\Gamma\subset\,G$. A lattice $\Gamma\subset\,G$ exists if and only if the Lie algebra $\lie{g}$ of $G$ has a rational structure, by a result of Mal'cev \cite{Malcev}. Nilmanifolds are parallelizable, hence, by the result of Crowley and Nordstr\"om mentioned in the introduction, they admit coclosed $\Gtwo$-structures. In this paper, however, we are interested in {\em left-invariant} (purely) coclosed $\Gtwo$-structures on nilmanifolds. By definition, these are left-invariant $\Gtwo$-structures on the corresponding Lie group.

\section{Constructing purely coclosed $\Gtwo$-structures}

In this section we review the construction of purely coclosed $\Gtwo$-structures on 7-dimen\-sional (nilpotent) Lie algebras presented in \cite{purely1}, to which we refer for the details and the proofs.

Let $\lie{g}$ be a $7$-dimensional Lie algebra with non-trivial center $\lie{z}(\fg)$ (recall that the center of a nilpotent Lie algebra is non-trivial). Let $V\subset\fg$ be a codimension 1 subspace, and let $X\in\lie{z}(\fg)$ be such that $\fg=V\oplus\langle X\rangle$. Let $\omega\in\Lambda^2\fg^*$ and $\psi_-\in\Lambda^3\fg^*$ be such that
\begin{itemize}
\item $\imath_X\omega=0$;
\item $\imath_X\psi_-=0$;
\item they define an $\SU(3)$-structure on $V$.
\end{itemize}

In particular, denoting by $\bar{\omega}$ and $\bar{\psi}_-$ the pull-back to $V$ of the above tensors, we have an almost symplectic 2-form $\bar{\omega}$ and a stable, negative 3-form $\bar{\psi}_-$ with $\bar{\omega}\wedge\bar{\psi}_-=0$. This determines $\bar{\psi}_+\in\Lambda^3(V^*)$. Extend $\bar{\psi}_+$ to an element $\psi_+\in\Lambda^3\fg^*$ by declaring $\imath_X\psi_+=0$. Finally, let $\eta\in\fg^*$ be such that
$\eta(X)\neq 0$.

It follows that $\varphi=\omega\wedge\eta+\psi_+$ is a $\Gtwo$-form on $\fg$; moreover, if $h$ denotes the induced $\SU(3)$-metric on $V$, then the $\Gtwo$-metric on $\fg$ is $g=g_\varphi=h+\eta\ox\eta$. Clearly, $*_\varphi\varphi=\frac{\omega^2}{2}+\psi_-\wedge\eta$. The following result, taken from \cite[Theorem 4.1]{purely1}, gives sufficient conditions on $(\omega,\psi_-,\eta)$ in order for the $\Gtwo$-structure $\varphi=\omega\wedge\eta+\psi_+$ to be (purely) coclosed.

\begin{theorem}\label{thm:construction}
In the above setting, the $\Gtwo$-structure is coclosed if
\begin{enumerate}
\item $d\psi_-=0$;
\item $\omega\wedge d\omega=\psi_-\wedge d\eta$.
\end{enumerate}
Furthermore, the coclosed $\Gtwo$-structure is pure if
\begin{enumerate}
\setcounter{enumi}{2}
\item $\omega^2 \wedge d\eta =-2 \psi_+ \wedge d\omega$.
\end{enumerate}
\end{theorem}

\begin{example}
Consider the Lie algebra $\mathfrak{g}=123457A=(0,0,-12,-13,-14,-15,-16)$ in the notation of \cite{Gong}. This means that $\mathfrak{g}$ is 7-dimensional and admits a basis $\{e_1,\ldots,e_7\}$, such that in terms of the dual basis $\{e^1,\ldots,e^7\}$, the Lie algebra structure is given by $de^1=de^2=0, \ de^3=-e^{12}, \ de^4=-e^{13}, \ de^5=-e^{14}, \ de^6=-e^{15}, \ de^7=-e^{16}$. Consider the subspace $V=\text{span}(e_1,e_2,e_3,e_4,e_5,e_6)\subset \mathfrak{g}$ and the central vector $X=e_7$. Set
\begin{itemize}
    \item $\omega=e^{12}+e^{34}-e^{56}$;
    \item $\psi_{-}=e^{135}+e^{236}+e^{146}-e^{245}$;
    \item $\eta=-e^3-e^7$.
\end{itemize}
Then $(\omega, \psi_{-})$ defines an $\text{SU}(3)$-structure on $V$ with induced metric $h=\sum_{i=1}^6 e^i \otimes e^i$, and 
$\psi_{+}=e^{235}+e^{145}-e^{136}+e^{246}$. Hence, $\phi=\omega \wedge \eta + \psi_{+}$ defines a $\Gtwo$-structure on $\mathfrak{g}$ with $\Gtwo$-metric $g=g_{\phi}=\sum_{i=1}^6 e^i\otimes e^i +e^3 \otimes e^3+e^3 \otimes e^7+ e^7 \otimes e^3+ e^7 \otimes e^7$.
Once checks that the conditions of Theorem \ref{thm:construction} are satisfied, the $\Gtwo$-structure is purely coclosed. 
\end{example}

In general, for the verification that $(\omega, \psi_{-})$ is an $\text{SU}(3)$-structure and that $\omega, \psi_{\pm}, \eta$ verify the condition of Theorem \ref{thm:construction} we refer to \cite{worksheets}. For an explanation of why the datasheets do compute the claimed quantities, we refer to \cite[Section $3$]{purely1}.

\section{Obstructions for coclosed $\Gtwo$-structures on nilpotent Lie algebras}\label{ostruzioni Giorgia}
In this section we present some obstructions to the existence of a coclosed $\Gtwo$-structure on a certain Lie algebra. As a result, we determine which 5-step nilpotent Lie algebras do not carry coclosed (hence purely coclosed) $\Gtwo$ sttructures. The first three obstructions were already introduced in \cite{purely1}, to which we refer for further details and for the proofs. The fourth obstruction is new.

We take the opportunity to correct a minor inaccuracy in the first obstruction presented in \cite{purely1}. Indeed, in \cite[Corollary 5.3]{purely1} one picks two vectors $X,Y\in\fg$ with $Y\in\mathfrak{z}(\fg)$ and considers the subspace $U\subset \Lambda^2\fg^*$ obtained by contracting a generic closed 4-form with $X$ and $Y$. The claim there is that if $\Lambda^2U=0$, then $\fg$ does not admit any coclosed $\Gtwo$-structure. In the proof it is claimed that if $\kappa\in\Lambda^4\fg^*$ is exact and $Y\in\mathfrak{z}(\fg)$ then $\iota_Y\kappa=0$. This is not true, in general. The first obstruction has to be stated in a slightly different way, which is actually a mere reformulation of \cite[Lemma 5.2]{purely1}.

\noindent\textbf{First obstruction} Let $\fg$ be a $7$-dimensional Lie algebra. Suppose that there exist linearly independent vectors $X,Y\in\fg$ such that, for every $\kappa\in\Lambda^4\fg^*$ closed, $\iota_X\iota_Y\kappa\in U$ for a subspace $U\subset\Lambda^2\fg^*$ such that $\Lambda^2 U=0$. Then, $\fg$ does not admit any coclosed $\Gtwo$-structure.

The two differences with respect the first obstruction in \cite{purely1} are:
\begin{itemize}
    \item one needs to check all closed forms, not just the representatives of cohomology classes;
    \item it is not necessary to pick one of the two vectors in the center of the Lie algebra.
\end{itemize}

We checked all the Lie algebras handled in \cite{purely1}, for which the existence of a coclosed $\Gtwo$-structure was excluded using the first obstruction. In all cases, we were able to find two vectors which satisfy the hypotheses of first obstruction above. Thus, the results of \cite{purely1} are correct.

The second obstruction is taken from \cite[Corollary 5.5]{purely1}.
 
\noindent\textbf{Second obstruction} 
Let $\lie{g}$ be a $7$-dimensional nilpotent Lie algebra and let $\{e_1,\ldots,e_7\}$ be a nilpotent basis, so that $[ e_i,e_j]=\sum_{k>i,j} c_{ij}^k e_k$ $\forall  i,j=1,\ldots,7$; let $\{e^1,\ldots,e^7\}$ be the dual basis. Take a list of generators of the space of closed $4$-forms $z_\alpha \in \Lambda^4(\lie{g}^*)$. Suppose that $z_\alpha \in \la e^1,e^2\ra \wedge \Lambda^3(\lie{g}^*)$ for all $\alpha$. Then $\lie{g}$ does not admit coclosed $\Gtwo$-structures.


The third obstruction appears in \cite[Proposition 5.6]{purely1}. Given a nilpotent basis $\{e_1,\ldots,e_7\}$ of $\fg$, fix $X=e_7$. We compute a basis $\{z_\alpha\}$ of the closed $3$-forms in $\lie{h}^*=\la e^1,\ldots,e^6\ra$. Therefore for a closed $3$-form, we can write $\tau= \sum a_\alpha z_\alpha$.

\noindent \textbf{Third obstruction}  Suppose we have elements $w_1,\ldots, w_\ell\in \Lambda^5(\lie{h}^*)$, and let $W$ be a subspace such that $\Lambda^5(\lie{h}^*)=W\oplus \la w_1,\ldots, w_\ell\ra$. Suppose furthermore that for every closed $2$-form $\beta$ and closed $3$-form $\tau$ on $\lie{h}^*$, we have 
\[
\beta\wedge d\beta ,\tau \wedge de^j \in W,  \ j=1,\ldots, 6\,.
\]
Then define the linear subspace
\[
H=\left\{ (a_\alpha) \, | \, \sum a_\alpha z_\alpha \wedge de^7 \in W\right\}.
\]
If $\lambda ( \sum a_\alpha z_\alpha)\geq 0$ for all $(a_\alpha)\in H$, then there is no coclosed $\Gtwo$-structure on $\lie{g}$; here $\lambda\colon\Lambda^3(\mathfrak{h}^*)\to (\Lambda^6(\mathfrak{h}^*))^{\otimes 2}$ is the so-called {\em Hitchin invariant}, introduced in \cite{Hitchin}.

\begin{remark}
    The obstructions have been checked using \texttt{SageMath}. For an explanation of why the ancillary data verify the second and third obstruction, we refer to \cite[Section 7]{purely1}. As for the ``new'' first obstruction, we need to contract all closed 4-forms with the two vectors $X$ and $Y$, not just the cohomology classes. For instance, for the Lie algebra 13457A we do:\end{remark}

{\tiny
\noindent \begin{verb}
A.<x1,x2,x3,x4,x5,x6,x7> = GradedCommutativeAlgebra(QQ)
\end{verb}

\noindent \begin{verb}
M=A.cdg_algebra({x3:-x1*x2, x4:-x1*x3, x5:-x1*x4, x7:-x1*x5-x2*x6})
\end{verb}

\noindent \begin{verb}
M.inject_variables()
\end{verb}

\noindent \begin{verb}
M.cohomology(4)
\end{verb}

\smallskip

\noindent \begin{verb}
Defining x1, x2, x3, x4, x5, x6, x7
\end{verb}

\smallskip

\noindent \begin{verb}
Free module generated by {[x1*x3*x4*x5], [x2*x3*x4*x5], [x2*x3*x5*x6 + x1*x3*x4*x7], [x1*x4*x5*x6], 
\end{verb}

\smallskip

\noindent \begin{verb}
[x2*x4*x5*x6 + x2*x3*x6*x7], [x3*x4*x5*x6 + x1*x4*x5*x7], [x1*x5*x6*x7]} over Rational Field
\end{verb}

\smallskip

\noindent \begin{verb}
D=[x1*x2*x3,x1*x2*x4,x1*x2*x5,x1*x2*x6,x1*x2*x7,x1*x3*x4,x1*x3*x5,x1*x3*x6,x1*x3*x7,x1*x4*x5,x1*x4*x6,
\end{verb}

\smallskip

\noindent \begin{verb}
x1*x4*x7,x1*x5*x6x1*x5*x7,x1*x6*x7,x2*x3*x4,x2*x3*x5,x2*x3*x6,x2*x3*x7,x2*x4*x5,x2*x4*x6,x2*x4*x7,
\end{verb}

\smallskip

\noindent \begin{verb}
x2*x5*x6,x2*x5*x7,x2*x6*x7,x3*x4*x5,x3*x4*x6,x3*x4*x7,x3*x5*x6,x3*x5*x7,x3*x6*x7,x4*x5*x6, x4*x5*x7,
\end{verb}

\smallskip

\noindent \begin{verb}
x4*x6*x7,x5*x6*x7]
\end{verb}

\smallskip

\noindent \begin{verb}
for i in range(35): print( "diff(",D[i] ,") =" ,D[i].differential())
\end{verb}

\smallskip

\noindent \begin{verb}
Defining x1, x2, x3, x4, x5, x6, x7
\end{verb}

\smallskip

\noindent \begin{verb}
diff( x1*x2*x3 ) = 0
\end{verb}

\smallskip

\noindent \begin{verb}
diff( x1*x2*x4 ) = 0
\end{verb}

\smallskip

\noindent \begin{verb}
diff( x1*x2*x5 ) = 0
\end{verb}

\smallskip

\noindent \begin{verb}
diff( x1*x2*x6 ) = 0
\end{verb}

\smallskip

\noindent \begin{verb}
diff( x1*x2*x7 ) = 0
\end{verb}

\smallskip

\noindent \begin{verb}
diff( x1*x3*x4 ) = 0
\end{verb}

\smallskip

\noindent \begin{verb}
diff( x1*x3*x5 ) = 0
\end{verb}

\smallskip

\noindent \begin{verb}
diff( x1*x3*x6 ) = 0
\end{verb}

\smallskip

\noindent \begin{verb}
diff( x1*x3*x7 ) = x1*x2*x3*x6
\end{verb}

\smallskip

\noindent \begin{verb}
diff( x1*x4*x5 ) = 0
\end{verb}

\smallskip

\noindent \begin{verb}
diff( x1*x4*x6 ) = 0
\end{verb}

\smallskip

\noindent \begin{verb}
diff( x1*x4*x7 ) = x1*x2*x4*x6
\end{verb}

\smallskip

\noindent \begin{verb}
diff( x1*x5*x6 ) = 0
\end{verb}

\smallskip

\noindent \begin{verb}
diff( x1*x5*x7 ) = x1*x2*x5*x6
\end{verb}

\smallskip

\noindent \begin{verb}
diff( x1*x6*x7 ) = 0
\end{verb}

\smallskip

\noindent \begin{verb}
diff( x2*x3*x4 ) = 0
\end{verb}

\smallskip

\noindent \begin{verb}
diff( x2*x3*x5 ) = -x1*x2*x3*x4
\end{verb}

\smallskip

\noindent \begin{verb}
diff( x2*x3*x6 ) = 0
\end{verb}

\smallskip

\noindent \begin{verb}
diff( x2*x3*x7 ) = -x1*x2*x3*x5
\end{verb}

\smallskip

\noindent \begin{verb}
diff( x2*x4*x5 ) = -x1*x2*x3*x5
\end{verb}

\smallskip

\noindent \begin{verb}
diff( x2*x4*x6 ) = -x1*x2*x3*x6
\end{verb}

\smallskip

\noindent \begin{verb}
diff( x2*x4*x7 ) = -x1*x2*x4*x5 - x1*x2*x3*x7
\end{verb}

\smallskip

\noindent \begin{verb}
diff( x2*x5*x6 ) = -x1*x2*x4*x6
\end{verb}

\smallskip

\noindent \begin{verb}
diff( x2*x5*x7 ) = -x1*x2*x4*x7
\end{verb}

\smallskip

\noindent \begin{verb}
diff( x2*x6*x7 ) = x1*x2*x5*x6
\end{verb}

\smallskip

\noindent \begin{verb}
diff( x3*x4*x5 ) = -x1*x2*x4*x5
\end{verb}

\smallskip

\noindent \begin{verb}
diff( x3*x4*x6 ) = -x1*x2*x4*x6
\end{verb}

\smallskip

\noindent \begin{verb}
diff( x3*x4*x7 ) = -x1*x3*x4*x5 - x2*x3*x4*x6 - x1*x2*x4*x7
\end{verb}

\smallskip

\noindent \begin{verb}
diff( x3*x5*x6 ) = -x1*x3*x4*x6 - x1*x2*x5*x6
\end{verb}

\smallskip

\noindent \begin{verb}
diff( x3*x5*x7 ) = -x2*x3*x5*x6 - x1*x3*x4*x7 - x1*x2*x5*x7
\end{verb}

\smallskip

\noindent \begin{verb}
diff( x3*x6*x7 ) = x1*x3*x5*x6 - x1*x2*x6*x7
\end{verb}

\smallskip

\noindent \begin{verb}
diff( x4*x5*x6 ) = -x1*x3*x5*x6
\end{verb}

\smallskip

\noindent \begin{verb}
diff( x4*x5*x7 ) = -x2*x4*x5*x6 - x1*x3*x5*x7
\end{verb}

\smallskip

\noindent \begin{verb}
diff( x4*x6*x7 ) = x1*x4*x5*x6 - x1*x3*x6*x7
\end{verb}

\smallskip

\noindent \begin{verb}
diff( x5*x6*x7 ) = -x1*x4*x6*x7
\end{verb}

\smallskip

\noindent \begin{verb}
# First obstruction: Use X=x5, Y=x7, U=<x1*x2,x1*x3,x1*x4,x1*x6>
\end{verb}}

\medskip

In the analysis carried out in this paper we bumped into two Lie algebras, $23457$B and $23457$F in Gong's notation, which passed the three obstructions above, but for which the method used in \cite{purely1} did not yield any result.

We thus develop a fourth obstruction, which allows us to prove that the above Lie algebras do not carry any coclosed $\Gtwo$-structure. This obstruction is based on the technique described in \cite{FinoMartinMRaffero} for computing the inner product determined by a positive $4$-form.

Let $\varphi$ be a $\Gtwo$-form, with associated metric $g_\varphi$, orientation $\mathrm{vol}_\varphi$ and Hodge-star $\ast_\varphi$. At each point of the manifold, $\varphi$ can be written as in \eqref{G2-form}. Recall that $\varphi$ is stable, in the sense of Hitchin (see \cite{Hitchin}). It turns out that the $4$-form $*_{\varphi}\varphi$ is also stable and positive; at each point, its $\text{GL}(7,\RR)$-orbit $\Lambda_{+}^4(\RR^7)^*$ is open in $\Lambda^4(\RR^7)^*$. The orbit $\Lambda_{+}^4(\RR^7)^*$ is isomorphic to
\[
\text{GL}(7,\RR) / \bigl\{\Gtwo \cup \bigl(\Gtwo \circ (-\mathrm{Id}_{\RR^7})\bigr)\bigr\},
\]
so the pointwise map
\[
\Lambda_{+}^3(\RR^7)^* \longrightarrow \Lambda_{+}^4(\RR^7)^*, \qquad
\varphi \longmapsto *_{\varphi}\varphi,
\]
is a double covering. Hence, a positive $4$-form defines an inner product but not an orientation. Consequently, in order to compute the inner product determined by a positive $4$-form, it is necessary to fix a background volume form. Once a volume form $\Omega$ on $\RR^7$ is fixed, there is an isomorphism
\[
\Lambda^4(\RR^7)^* \cong \Lambda^3(\RR^7) \otimes \Lambda^7(\RR^7)^*
\]
that identifies $\rho \in \Lambda^4(\RR^7)^*$ with the element
$\hat{\rho}\otimes \Omega \in \Lambda^3(\RR^7)\otimes \Lambda^7(\RR^7)^*$ such that 
\begin{equation}\label{step:2}
 \iota_{\hat{\rho}}\Omega = \rho\,.
\end{equation}
One then obtains a symmetric bilinear map
\begin{align}\label{metricmatrix}
    \mathrm{B}_{\rho}: (\RR^7)^* \times (\RR^7)^*
    &\longrightarrow \Lambda^7(\RR^7)\otimes \bigl(\Lambda^7(\RR^7)^*\bigr)^{\otimes 3}
       \cong \bigl(\Lambda^7(\RR^7)^*\bigr)^{\otimes 2} \nonumber \\
    (v,w) &\longmapsto
    \left( \frac{1}{6}\,\iota_v \hat{\rho}\wedge \iota_w \hat{\rho} \wedge \hat{\rho}\right)
    \otimes \Omega^{\otimes 3}.
\end{align}
The $4$-form $\rho$ is stable if and only if
\[
\bigl|\det(\mathrm{B}_{\rho})\bigr|^{\frac{1}{12}} \in \Lambda^7(\RR^7)^*
\]
is non-zero. 


\begin{remark}
    The $\Gtwo$-metric induced by a definite $4$-form $\rho$ is then defined as $g_\rho=|\det(\mathrm{B}_{\rho})|^{-\frac{1}{6}}\mathrm{B}_{\rho}$. As noticed in \cite{FinoMartinMRaffero}, if $\varphi$ is a $\Gtwo$-form, the relation between $g_\varphi$ and the metric $g_{\ast_\varphi\varphi}$ induced by the $4$-form $\ast_\varphi\varphi$ is given by $g_{\ast_\varphi\varphi}= g_\varphi^{-1}$.
\end{remark}

The Lie algebras we deal with are described in terms of a (nilpotent) basis $\{e_1,\ldots,e_7\}$, hence $\Omega=e^{1234567}$ is a natural choice for the background volume form. We can state the fourth obstruction:

\begin{proposition}[Fourth obstruction]\label{4th obs Giorgia}
Let $\mathfrak{g}$ be a 7-dimensional Lie algebra, and let $\{e_1,\ldots,e_7\}$ be a basis of $\fg$. If, for every closed $4$-form $\rho\in\Lambda^4\fg^*$, the matrix $B_{\rho}$ 
defined in \eqref{metricmatrix} fails to be positive definite, then $\fg$ admits no coclosed $\mathrm{G}_2$-structure.
\end{proposition}

The computations needed to construct the matrix $\mathrm{B}_{\rho}$ are implemented in \texttt{SageMath} and can be found in \cite{worksheets}, to which we refer the reader. Here is a summary of the steps we perform in order to check the fourth obstruction:
\begin{enumerate}
    \item The first step is find the closed 4-forms in the Lie algebra $\fg$ under examination. We compute separately the cohomology classes in $H^4(\mathfrak{g})$ (step 1a) and the exact 4-forms (step 1b).
    \item We fix the background volume form $\Omega$ using the nilpotent basis in which $\fg$ is given and compute $\hat{\rho}$ using \eqref{step:2}.
    \item Once we have $\hat{\rho}$ we could compute $\textrm{B}_\rho$ using \eqref{metricmatrix}. However, we want to show that no choice of a closed 4-form $\rho$ produces a positive-definite $\textrm{B}_\rho$. To do this, we use Sylvester's criterion and compute the $(1,1)$ entry of $\textrm{B}_\rho$, which turns out to be 0 in both cases. This of course prevents the existence of a stable 4-form.
\end{enumerate}

\begin{theorem}\label{5step-indecomposable-non-coclosed}
The following indecomposable 5-step nilpotent Lie algebras do not admit any coclosed $\Gtwo$-structure: $12457$F, $12457$L, $13457$A, $13457$B, $13457$C, $13457$E, $13457$G, $13457$I, $23457$A, $23457$B, $23457$F.
\end{theorem}
\begin{proof}
    The proof is based on the obstructions mentioned above; the relevant computations can be found in the ancillary worksheets \cite{worksheets}. For the convenience of the reader, Table \ref{table obs 5-sp Gior} contains the obstruction used in each case.
    \begin{table}[H]
\centering
\caption{Obstructions used to rule out the existence of coclosed $\Gtwo$-structures on indecomposable 5-step NLAs}\label{table obs 5-sp Gior}
\vspace{0.25 cm}
{\tabulinesep=1.2mm
\begin{tabu}{l|c||l|c}
\toprule[1.5pt]
NLA & Obstruction & NLA & Obstruction\\

\specialrule{1pt}{0pt}{0pt}
\text{$12457$F} & Third & \text{$13457$G}& Third\\
\specialrule{1pt}{0pt}{0pt}
\text{$12457$L} & Third & \text{$13457$I}& Third\\
\specialrule{1pt}{0pt}{0pt}
\text{$13457$A}& First & \text{$23457$A}& Second\\
\specialrule{1pt}{0pt}{0pt}
\text{$13457$B}& First & \text{$23457$B}& Fourth\\
\specialrule{1pt}{0pt}{0pt}
\text{$13457$C}& First & \text{$23457$F}& Fourth\\
\specialrule{1pt}{0pt}{0pt}
\text{$13457$E}& Third\\
\bottomrule[1pt]
\end{tabu}}
\end{table} 
\end{proof}

\section{Seven-dimensional nilpotent lie algebras with purely coclosed $\Gtwo$-structures}
\subsection{Indecomposable 5-step nilpotent Lie algebras}
The seven-dimensional indecomposable 5-step nilpotent Lie algebras are listed in Table \ref{7d-ind-5step-1}. According to Gong’s classification \cite{Gong}, 
there are 36 of them. Two of them, $12457N(\lambda)$ and $12457N_2(\lambda)$, depend on a real parameter; the conditions on  $\lambda$ are given in Table \ref{7d-ind-5step-1}.

 \begin{theorem}\label{5step-indecomposable-pure}
All $7$-dimensional indecomposable 5-step nilpotent Lie algebras admit a purely coclosed $\Gtwo$-structure, except for: $12457$F, $12457$L, $13457$A, $13457$B, $13457$C, $13457$E, $13457$G, $13457$I, $23457$A, $23457$B, $23457$F.
\end{theorem}

\begin{proof}
    The fact that the listed Lie algebras do not admit purely coclosed $\Gtwo$ structures was proved in the above Theorem \ref{5step-indecomposable-non-coclosed}. Tables \ref{table:Giorg1} and \ref{table:Giorg2} contain the purely coclosed $\Gtwo$-structures that we have found on the remaining 5-step indecomposable NLA; all of them are constructed according to Theorem \ref{thm:construction}.
\end{proof}

\subsection{Indecomposable 6-step nilpotent Lie algebras}
We deal with seven-dimensional indecomposable 6-step nilpotent Lie algebras. According to Gong \cite{Gong}, there are 9  of them, listed in Table \ref{7d-ind-6step-1}. Only the algebra $123457\text{I}(\lambda)$ depends on a real parameter $\lambda$. In the following theorem, we present the results we obtained.

\begin{theorem}\label{6step-indecomposable-pure}
All $7$-dimensional (indecomposable) 6-step nilpotent Lie algebras admit a purely coclosed $\Gtwo$-structure.
\end{theorem}

The purely coclosed $\Gtwo$-structures that we have constructed on these Lie algebras are explicitly given in Table \ref{table:Giorg3}.

\section{Further comments}

Using results from \cite{BMR}, it was proved in \cite[Theorem 6.2]{purely1} that every seven-dimensional decomposable nilpotent Lie algebra admitting a coclosed $\Gtwo$-structure also admits a purely coclosed one, except for $\fh_3\oplus\R^4$, the direct sum of the 3-dimensional Heisenberg algebra and $\R^4$ (this is denoted $\fn_2$ in \cite{purely1}). In the indecomposable case, it was observed in \cite[Corollary 6.8]{purely1} that every seven-dimensional indecomposable nilpotent Lie algebra of nilpotency step $\leq 4$ admitting a coclosed $\Gtwo$-structure also admits a purely coclosed one. 

As a consequence of the results presented in the previous section, we can complete the picture in the indecomposable case, obtaining the following theorem.

\begin{theorem}\label{thm:nice}
Let $\fg$ be a 7-dimensional nilpotent Lie algebra.
\begin{itemize}
    \item If $\fg$ is decomposable, then it admits a coclosed $\Gtwo$-structure if and only if it admits a purely coclosed one, unless $\fg=\fh_3\oplus\R^4$.
    \item If $\fg$ is indecomposable, then it admits a coclosed $\Gtwo$-structure if and only if it admits a purely coclosed one.
\end{itemize}
\end{theorem}

\begin{remark}
It seems difficult to deduce the existence of a purely coclosed $\Gtwo$-structure directly from the existence of a coclosed one, under the indecomposability assumption. Also, it seems hard to deduce directly the exceptionality of $\fh_3\oplus\R^4$ in this context.
\end{remark}

\begin{remark}
It is worth pointing out that all purely coclosed $\Gtwo$-structures have been constructed in a uniform way, using the {\em Ansatz} of Theorem \ref{thm:construction}.
\end{remark}

The class of coclosed $\Gtwo$-structure contains both the pure classes of nearly parallel and purely coclosed $\Gtwo$-structures. Clearly, the $\Gtwo$ 4-form is exact in the nearly parallel case. Using \texttt{SageMath} we checked that none of the purely coclosed $\Gtwo$ 4-forms we constructed is exact. It is an open question whether the $\Gtwo$ 4-form of a purely coclosed $\Gtwo$-structure can be exact. A similar question exists in the context of closed $\Gtwo$ structures: it is not known whether the $\Gtwo$ 3-form of a closed $\Gtwo$-structure can be exact. Negative results on compact quotients of Lie groups have been obtained in \cite{FMMR}.

In \cite{LUV}, the authors gave a positive answer as to whether there exist infinitely many non-isomorphic nilpotent Lie algebras in dimension 8 admitting complex structures. Recall that a nilpotent Lie algebra admits a K\"ahler structure if and only if it is abelian. Now in dimension 7, a nilpotent Lie algebra can not admit torsion-free, nearly parallel nor locally conformally parallel $\Gtwo$-structures. Thus, a way to transfer the above question to the 7-dimensional world could be the following:
\begin{question}
Do there exist infinitely many non-isomorphic nilpotent Lie algebras in dimension 7 admitting closed or (purely) coclosed $\Gtwo$-structures?
\end{question}

Using results of \cite{purely1,CF} and from the above analysis, we are able to give an answer to this question:
\begin{itemize}
    \item in the closed setting, there are only finitely many non-isomorphic nilpotent Lie algebras which admit closed $\Gtwo$-structures; in particular, there exist only finitely many {\em real homotopy types} of 7-dimensional nilmanifolds admitting a closed $\Gtwo$-structure;
    \item in the purely coclosed setting, for every nilpotency step $k\geq 3$ there are infinitely many non-isomorphic nilpotent Lie algebras which admit purely closed $\Gtwo$-structures. Explicit examples are given by the following families:
    \begin{itemize}
        \item $147\textrm{E}(\lambda)$, with nilpotency step 3;
        \item $1357\textrm{M}(\lambda)$, with nilpotency step 4;
        \item $12457\textrm{N}(\lambda)$, with nilpotency step 5;
        \item $123457\textrm{I}(\lambda)$, with nilpotency step 6.
    \end{itemize} 
    In particular, there exist infinitely many {\em real homotopy types} of 7-dimensional nilmanifolds admitting a purely coclosed $\Gtwo$-structure.
\end{itemize}

It not known to the authors whether there exist infinitely many {\em rational} homotopy types of nilmanifolds admitting a closed $\Gtwo$-structure.

\newpage

\begin{table}[H]
\centering
\caption{Purely coclosed $\Gtwo$-structures on indecomposable 5-step NLAs $-$ I}\label{table:Giorg1}
\vspace{0.25 cm}
{\tabulinesep=1.2mm
\begin{tabu}{l|c|c|c}
\toprule[1.5pt]
NLA & $\omega$ & $\psi_-$ & $\eta$\\
\specialrule{1pt}{0pt}{0pt}
\text{$12357$A}
    &$\begin{array}{c}-e^{13}+e^{24}-\frac{2}{5}e^{56}\\[1.5pt]+\frac{3}{5}e^{45}+\frac{3}{5}e^{26}\end{array}$& $e^{125}+e^{236}+e^{146}+e^{345}$ & $e^7$\\
\specialrule{1pt}{0pt}{0pt}
\text{$12357$B} &$\begin{array}{c}-e^{13}+e^{24}-\frac{2}{5}e^{56}\\[1.5pt]+\frac{3}{5}e^{45}+\frac{3}{5}e^{26}\end{array}$& $e^{125}+e^{236}+e^{146}+e^{345}$ & $-e^{5}+2e^7$\\
\specialrule{1pt}{0pt}{0pt}
\text{$12357\text{B}_1$}
    &$\begin{array}{c}-e^{13}+e^{24}-\frac{2}{5}e^{56}\\[1.5pt]+\frac{3}{5}e^{45}+\frac{3}{5}e^{26}\end{array}$& $e^{125}+e^{236}+e^{146}+e^{345}$ & $e^{5}+2e^7$\\
\specialrule{1pt}{0pt}{0pt}
\text{$12357$C} & $-e^{13}+e^{24}-e^{56}$& $e^{125}+e^{236}+e^{146}+e^{345}$ &$-3e^7 $\\
\specialrule{1pt}{0pt}{0pt}
\text{$12457$A}& $-e^{15}-e^{24}-e^{36}+e^{35}$
& \text{$e^{123}+e^{145}-e^{146}+2e^{256}+e^{345}$} &$e^7$\\
\specialrule{1pt}{0pt}{0pt}
\text{$12457$B}& $-e^{15}-e^{24}-e^{36}+e^{35}$ 
& \text{$e^{123}+e^{145}-e^{146}+2e^{256}+e^{345}$} &$e^6+e^7$\\
\specialrule{1pt}{0pt}{0pt}
\text{$12457$C} & $\begin{array}{c}e^{14}+e^{15}+e^{24}\\[1.5pt] +\frac{1}{2}e^{25}-e^{35}-e^{36}\end{array}$ 
& \text{$e^{123}+2e^{145}+2e^{146}-e^{256}+e^{345}$} & $-\frac{3}{2}e^3-\frac{1}{2}e^6-\frac{1}{4}e^7$\\
\specialrule{1pt}{0pt}{0pt}
\text{$12457$D} & $e^{15}+e^{24}-e^{36}$ & \text{$e^{123}+2e^{146}-e^{256}+e^{345}$} & $e^7$\\
\specialrule{1pt}{0pt}{0pt}
\text{$12457$E} & $-e^{15}+e^{23}-e^{24}-e^{36}$ & \text{$\begin{array}{c}e^{123}+e^{124}+e^{125}-e^{146}\\[1.5pt]+e^{235}+e^{236}+2e^{256}+e^{345}\end{array}$} & $\frac{4}{3}e^3-\frac{2}{3}e^4+\frac{4}{3}e^7$\\
\specialrule{1pt}{0pt}{0pt} \text{$12457$G} & $e^{15}+e^{24}-e^{36}+2e^{13}$ & \text{$e^{123}+2e^{146}-e^{256}+e^{345}-2e^{236}$} & $2e^4+e^7$\\
\specialrule{1pt}{0pt}{0pt}
\text{$12457$H} & $e^{12}+e^{34}-e^{56}+e^{16}$ & \text{$-e^{135}-2e^{146}-e^{236}+e^{245}-e^{124}$} & $e^3+e^7$\\
\specialrule{1pt}{0pt}{0pt}
\text{$12457$I} & $e^{12}+e^{34}-e^{56}$ & \text{$-e^{135}-2e^{146}-e^{236}+e^{245}$} & $3e^3-e^7$\\
\hline \text{$12457$J} & $e^{12}-e^{34}+e^{56}-e^{14}$ & \text{$e^{135}+e^{146}-e^{236}+e^{245}+e^{126}$} & $-\frac{1}{2}e^3-e^5+\frac{1}{2}e^7$\\
\specialrule{1pt}{0pt}{0pt}
\text{$12457\text{J}_1$} & $e^{12}+e^{34}-e^{56}$ & \text{$e^{135}+e^{146}+e^{236}-e^{245}$} & $-\frac{3}{2}e^3+\frac{1}{2}e^5-\frac{1}{2}e^7$\\
\specialrule{1pt}{0pt}{0pt}
\text{$12457$K} & $e^{12}-e^{34}+e^{56}-e^{14}$ & \text{$e^{135}+e^{146}-e^{236}+e^{245}+e^{126}$} & $-e^5+e^7$\\
\specialrule{1pt}{0pt}{0pt} \text{$12457\text{L}_1$} & 
  $\begin{array}{c}e^{12}+e^{23}+e^{24}+2e^{35}\\[1.5pt]+2e^{45}+e^{56}+e^{46}
  \end{array}$  & $\begin{array}{c}e^{136}+2e^{234}-e^{145}-2e^{256}\\[1.5pt]-e^{236}+e^{345}-e^{245}+e^{146}\end{array}$ 
& $\begin{array}{c}
-\frac{2}{3}e^4-\frac{5}{3}e^5\\[1.5pt]+\frac{7}{3}e^6+\frac{1}{3}e^7\end{array}$\\
\specialrule{1pt}{0pt}{0pt}
$\begin{array}{c}12457\text{N}(\lambda)\\[1.5pt]\lambda \neq \frac{78}{331}\end{array}$ & $\begin{array}{c}-e^{12}+2e^{14}\\[1.5pt]-6e^{15}-3e^{34}\\[1.5pt]-e^{36}-9e^{46}-e^{56}\end{array}$ &$ \begin{array}{c} -e^{235}-3e^{145}-3e^{136}+e^{246}\\[1.5pt]+e^{156}-2e^{345}-3e^{146}-3e^{256}\end{array}$ &$\begin{array}{c}\frac{2}{7}\frac{3570\lambda - 5329}{331\lambda -78}e^3\\[1.5pt]-\frac{1}{7}\frac{40105 \lambda +10926}{331\lambda -78}e^4\\[1.5pt]-\frac{1}{7}\frac{126\lambda + 7369}{331\lambda -78}e^5\\[1.5pt]-\frac{117}{7}\frac{485\lambda+7}{331\lambda -78}e^6\\[1.5pt]+\frac{1}{7}\frac{40147}{331\lambda -78}e^7\end{array}$\\
\bottomrule[1pt]
\end{tabu}}
\end{table} 
 \newpage

\begin{table}[H]
\centering
\caption{Purely coclosed $\Gtwo$-structures on indecomposable 5-step NLAs $-$ II}\label{table:Giorg2}
\vspace{0.25 cm}
{\tabulinesep=1.2mm
\begin{tabu}{l|c|c|c}
\toprule[1.5pt]
NLA & $\omega$ & $\psi_-$ & $\eta$\\
\specialrule{1pt}{0pt}{0pt}
$12457\text{N}\left(\frac{78}{331}\right)$ 
&$\begin{array}{c}
-e^{12}+4e^{14}\\[1.5pt]-12e^{15}-3e^{34}\\[1.5pt]-2e^{36}-18e^{46}-e^{56}\end{array}$
&$\begin{array}{c}-e^{235}-3e^{145}-3e^{136}+2e^{246}\\[1.5pt]+2e^{156}-4e^{345}-6e^{146}-6e^{256}\end{array}$
     &$\begin{array}{c}-\frac{55669637}{6091905}e^3\\[1.5pt]-\frac{2363031}{86410}e^4\\[1.5pt]-\frac{118524151}{6091905}e^5\\[1.5pt]-\frac{264950153}{2030635}e^6\\[1.5pt]+\frac{724456721}{6091905}e^7 \end{array}$\\
\specialrule{1pt}{0pt}{0pt}
\text{$12457\text{N}_1$}
    &$\begin{array}{c}e^{12}+e^{23}\\[1.5pt]+e^{24}+2e^{35}\\[1.5pt]+2e^{45}+e^{56}+e^{46}\end{array}$
        &$\begin{array}{c}e^{136}+2e^{234}-e^{145}-2e^{256}\\[1.5pt]-e^{236}+e^{345}-e^{245}+e^{146}\end{array}$  & $\begin{array}{c}-\frac{1}{3}e^4-\frac{4}{3}e^5\\[1.5pt]+2e^6+\frac{1}{3}e^7\end{array}$\\
\specialrule{1pt}{0pt}{0pt}
\text{$12457\text{N}_2(\lambda)$} &$\begin{array}{cc}3e^{12}+e^{13}\\[1.5pt]+2e^{15}+2e^{35}\\[1.5pt]-e^{36}+e^{45}+e^{46}\end{array}$
&$\begin{array}{c}e^{136}-e^{145}+e^{146}+2e^{234}\\[1.5pt]-e^{236}-e^{245}-2e^{256}+e^{345}\end{array}$ &
$\begin{array}{c}\frac{1}{32}(-12e^3\\[1.5pt]+(17\lambda-48)e^4\\[1.5pt]+(17\lambda-126)e^5\\[1.5pt]+(-17\lambda+113)e^6\\[1.5pt]+17e^7)\end{array}$\\
\specialrule{1pt}{0pt}{0pt}
\text{$13457$D}& $-e^{12}+e^{34}-e^{36}+e^{56}$& $-e^{136}-e^{145}-e^{156}+e^{235}-e^{246}$ &$\frac{1}{3}e^4-\frac{1}{3}e^7$\\
\specialrule{1pt}{0pt}{0pt}
\text{$13457$F} & $-e^{12}+e^{35}-e^{46}$ & \text{$\begin{array}{c}
e^{134}+e^{136}+e^{145}\\[1.5pt]+e^{156}-e^{236}-e^{245}\end{array}$} & $\begin{array}{c}
-2e^3-e^4+e^5\\[1.5pt] +e^6+e^7\end{array}$\\ 
\specialrule{1pt}{0pt}{0pt}
\text{$23457$C} & $-e^{16}+e^{25}+e^{34}-2e^{56}$ & $-e^{124}+e^{135}+e^{146}+e^{236}-e^{245}$ & $-e^3+2e^7$\\ 
\specialrule{1pt}{0pt}{0pt}
\text{$23457$D} & $e^{12}+e^{34}-e^{56}$ & 
$e^{135}+e^{146}-e^{235}+e^{236}-e^{245}$ & $e^3+e^7$\\ 
\specialrule{1pt}{0pt}{0pt}
\text{$23457$E} & $-e^{12}+e^{34}-e^{45}+e^{57}$ 
& \text{$\begin{array}{c}
e^{137}-e^{235}+e^{145}\\[1.5pt]+e^{247}+e^{157}-e^{245}\end{array}$} & $\begin{array}{c}-\frac{8}{3}e^3-\frac{1}{2}e^4-\frac{1}{2}e^5\\[1.5pt]+\frac{1}{2}e^6+\frac{1}{2}e^7\end{array}$ \\ 
\specialrule{1pt}{0pt}{0pt}
\text{$23457$G} & $\begin{array}{c} -e^{15}-e^{24}+2e^{34}\\[1.5pt]-e^{36}-e^{45}-2e^{56}\end{array}$ & $\begin{array}{c}-e^{123}-e^{125}+e^{135}\\[1.5pt]+e^{146}+e^{236}-e^{245}\end{array}$ & $\begin{array}{c}\frac{2}{3}e^3-3e^4\\[1.5pt]+2e^6+3e^7\end{array}$\\ 
\bottomrule[1pt]
\end{tabu}}
\end{table}

\begin{table}[H]
\centering
\caption{Purely coclosed $\Gtwo$-structures on indecomposable 6-step NLAs}\label{table:Giorg3}
\vspace{0.25 cm}
{\tabulinesep=1.2mm
\begin{tabu}{l|c|c|c}
\toprule[1.5pt]
NLA & $\omega$ & $\psi_-$ & $\eta$\\
\specialrule{1pt}{0pt}{0pt} \text{$123457$A}&$e^{12}+e^{34}-e^{56}$ & $e^{135}+e^{236}+e^{146}-e^{245}$ &$-e^3-e^7$\\
\specialrule{1pt}{0pt}{0pt} \text{$123457$B}&$e^{12}+e^{34}-e^{56}$ & $e^{135}+e^{236}+e^{146}-e^{245}$ &$-e^3+e^5-e^7$\\
\specialrule{1pt}{0pt}{0pt} \text{$123457$C}&$e^{12}+e^{34}-e^{56}$ & $e^{135}+e^{236}+2e^{146}-e^{245}$ &$\frac{1}{2}e^3+e^7$\\
\specialrule{1pt}{0pt}{0pt} \text{$123457$D}&$e^{12}+e^{34}-e^{56}$ & $e^{135}+e^{236}+e^{146}-e^{245}-e^{235}$ &$\begin{array}{c}-2e^4-e^5\\[1.5pt]+e^6-e^7\end{array}$\\
\specialrule{1pt}{0pt}{0pt} \text{$123457$E}&$e^{12}+e^{34}-e^{56}$ & $e^{135}+e^{236}+e^{146}-e^{245}-e^{235}$ &$-e^4+e^6-e^7$\\
\specialrule{1pt}{0pt}{0pt} \text{$123457$F}&$e^{12}+e^{34}-e^{56}$ & $e^{135}+2e^{236}+e^{146}-2e^{245}-e^{235}$ &$\begin{array}{c}-\frac{9}{5}e^3-\frac{5}{8}e^4\\[1.5pt]-\frac{1}{4}e^5+\frac{1}{2}e^6-e^7\end{array}$\\
\specialrule{1pt}{0pt}{0pt} \text{$123457$H}&$e^{12}+e^{34}-e^{56}$ & $e^{135}+e^{236}+2e^{146}-e^{245}$ &$-\frac{5}{2}e^3-\frac{2}{3}e^5+e^7$\\
\specialrule{1pt}{0pt}{0pt} \text{$123457\text{H}_1$}&$e^{12}+e^{34}-e^{56}$ & $e^{135}+e^{236}+2e^{146}-e^{245}$ &$-\frac{5}{2}e^3+\frac{2}{3}e^5-e^7$\\
\specialrule{1pt}{0pt}{0pt} \text{$123457\text{I}(\lambda)$, $\lambda \neq \frac{1}{2}$}&$e^{12}+e^{34}-e^{56}$ & $e^{135}+e^{236}+2e^{146}-e^{245}$ &$\begin{array}{c}\frac{1}{2}(\frac{2(\lambda -1)}{2 \lambda -1}-5)e^3\\[1.5pt]+\frac{1}{2\lambda-1}e^7\end{array}$\\
\specialrule{1pt}{0pt}{0pt} \text{$123457\text{I}(\lambda)$, $\lambda= \frac{1}{2}$}&$e^{12}+e^{34}-e^{56}$ & $e^{135}+e^{236}+3e^{146}-e^{245}$ &$-\frac{10}{3}e^3+2e^7$\\
\bottomrule[1pt]
\end{tabu}}
\end{table}


\appendix

\section{7D nilpotent Lie alegbras of nilpotency step $\ge 5$}\label{appendix}
In this appendix we list the seven-dimensional nilpotent Lie algebras of nilpotency steps $5$ and $6$. In Gong’s classification \cite{Gong}, the structure equations are expressed in terms of the Lie bracket $[ \cdot, \cdot]\colon\Lambda^2 \fg \to \fg$ with respect to a nilpotent basis $\{e_1,\ldots,e_7\}$. In the lists below, the structure equations are instead presented in terms of the Chevalley–Eilenberg differential $d\colon\fg^* \to \Lambda^2 \fg^*$, which is the dual operator of the Lie bracket, with respect to the dual nilpotent basis $\{e^1,\ldots,e^7\}$ of $\mathfrak{g}^*$. The notation $(0^3,-12,-14-23,-15+34,-16+35)$ means that, with respect to the ordered basis $\{e_1,\ldots,e_7\}$, the differential is given as follows: $de^i=0$ for $i=1,2,3$ and
\[
de^4=-e^{12}, \quad de^5=-e^{14}-e^{23}, \quad de^6=-e^{15}+e^{34}, \quad de^7=-e^{16}+e^{35}\,,
\]
where $e^{ij}$ is a short-hand notation for $e^i\wedge e^j$ and similarly for $e^{ijk\ell\cdots}$.

{\footnotesize 
\begin{table}[H]
\centering
\caption{7-dimensional indecomposable 5-step nilpotent Lie algebras}\label{7d-ind-5step-1}
\vspace{0.25 cm}
{\tabulinesep=1.2mm
\begin{tabu}{l|c|c}
\toprule[1.5pt]
NLA & structure equations & center\\
\specialrule{1pt}{0pt}{0pt}
\text{$12357$A} & $(0^3,-12,-14-23,-15+34,-16+35)$ & $\langle e_7 \rangle$\\
\specialrule{1pt}{0pt}{0pt}
\text{$12357$B} & $(0^3,-12,-14-23,-15+34,-16-23+35)$ & $\langle e_7 \rangle$\\
\specialrule{1pt}{0pt}{0pt} 
\text{$12357\text{B}_1$} & $(0^3,-12,-14-23,-15+34,-16+23+35)$ & $\langle e_7 \rangle$\\
\specialrule{1pt}{0pt}{0pt}
\text{$12357$C} & $(0^3,-12,-14-23,-15+34,-16-24+35)$ & $\langle e_7 \rangle$\\
\specialrule{1pt}{0pt}{0pt}
\text{$12457$A} & $(0^2,-12,-13,0,-14-25,-16-35)$ & $\langle e_7 \rangle$\\
\specialrule{1pt}{0pt}{0pt}
\text{$12457$B} & $(0^2,-12,-13,0,-14-25,-16-25-35)$ & $\langle e_7 \rangle$\\
\specialrule{1pt}{0pt}{0pt}
\text{$12457$C} & $(0^2,-12,-13,0,-14-25,-26+34)$ & $\langle e_7 \rangle$\\
\specialrule{1pt}{0pt}{0pt}
\text{$12457$D} & $(0^2,-12,-13,0,-14-25,-15-26+34)$ & $\langle e_7 \rangle$\\
\specialrule{1pt}{0pt}{0pt}
\text{$12457$E} & $(0^2,-12,-13,0,-14-23-25,-16-24-35)$ & $\langle e_7 \rangle$\\
\specialrule{1pt}{0pt}{0pt}
\text{$12457$F} & $(0^2,-12,-13,0,-14-23-25,-26+34)$ & $\langle e_7 \rangle$\\
\specialrule{1pt}{0pt}{0pt}
\text{$12457$G} & $(0^2,-12,-13,0,-14-23-25,-15-26+34)$ & $\langle e_7 \rangle$\\
\specialrule{1pt}{0pt}{0pt}
\text{$12457$H} & $(0^2,-12,-13,-23,-15-24,-16-34)$ & $\langle e_7 \rangle$\\
\specialrule{1pt}{0pt}{0pt}
\text{$12457$I} & $(0^2,-12,-13,-23,-15-24,-16-25-34)$ & $\langle e_7 \rangle$\\
\specialrule{1pt}{0pt}{0pt}
\text{$12457$J} & $(0^2,-12,-13,-23,-15-24,-14-16-25-34)$ & $\langle e_7 \rangle$\\
\specialrule{1pt}{0pt}{0pt}
\text{$12457\text{J}_1$} & $(0^2,-12,-13,-23,-15-24,-14-16+25-34)$ & $\langle e_7 \rangle$\\
\specialrule{1pt}{0pt}{0pt}
\text{$12457$K} & $(0^2,-12,-13,-23,-15-24,-14-16-34)$ & $\langle e_7 \rangle$\\
\specialrule{1pt}{0pt}{0pt}
\text{$12457$L} & $(0^2,-12,-13,-23,-15-24,-16-26-34+35)$ & $\langle e_7 \rangle$\\
\specialrule{1pt}{0pt}{0pt}
\text{$12457\text{L}_1$} & $(0^2,-12,-13,-23,14+25,-16+35)$ & $\langle e_7 \rangle$\\
\specialrule{1pt}{0pt}{0pt}
\text{$12457\text{N}(\lambda)$} & $(0^2,-12,-13,-23,-15-24,-14-16-\lambda\cdot 25-26-34+35)$ & $\langle e_7 \rangle$\\
\specialrule{1pt}{0pt}{0pt}
\text{$12457\text{N}_1$} & $(0^2,-12,-13,-23,14+25,-16-25+35)$ & $\langle e_7 \rangle$\\
\specialrule{1pt}{0pt}{0pt}
\text{$12457\text{N}_2(\lambda)$ $\lambda \geq 0$} & $(0^2,-12,-13,-23,14+25,-15-16-24-\lambda \cdot 25 +35)$ & $\langle e_7 \rangle$\\
\specialrule{1pt}{0pt}{0pt} 
\text{$13457$A} & $(0^2,-12,-13,-14,0,-15-26)$ & $\langle e_7 \rangle$\\
\specialrule{1pt}{0pt}{0pt}
\text{$13457$B} & $(0^2,-12,-13,-14,0,-15-23-26)$ & $\langle e_7 \rangle$\\
\specialrule{1pt}{0pt}{0pt}
\text{$13457$C} & $(0^2,-12,-13,-14,0,-16-25+34)$ & $\langle e_7 \rangle$\\
\specialrule{1pt}{0pt}{0pt}
\text{$13457$D} & $(0^2,-12,-13,-14-23,0,-15-24-26)$ & $\langle e_7 \rangle$\\
\specialrule{1pt}{0pt}{0pt}
\text{$13457$E} & $(0^2,-12,-13,-14-23,0,-16-25+34)$ & $\langle e_7 \rangle$\\
\specialrule{1pt}{0pt}{0pt}
\text{$13457$F} & $(0^2,-12,-13,-14,-23,-15-26)$ & $\langle e_7 \rangle$\\
\specialrule{1pt}{0pt}{0pt}
\text{$13457$G} & $(0^2,-12,-13,-14,-23,-16-24-25+34)$ & $\langle e_7 \rangle$\\
\specialrule{1pt}{0pt}{0pt}
\text{$13457$I} & $(0^2,-12,-13,-14,-23,-15-25-26+34)$ & $\langle e_7 \rangle$\\
\specialrule{1pt}{0pt}{0pt}
\text{$23457$A} & $(0^2,-12,-13,-14,-15,-23)$ & $\langle e_6, e_7 \rangle$\\
\specialrule{1pt}{0pt}{0pt}
\text{$23457$B} & $(0^2,-12,-13,-14,-25+34,-23)$ & $\langle e_6, e_7 \rangle$\\
\specialrule{1pt}{0pt}{0pt}
\text{$23457$C} & $(0^2,-12,-13,-14,-15,-25+34)$ & $\langle e_6,e_7 \rangle$ \\ 
\specialrule{1pt}{0pt}{0pt}
\text{$23457$D} & $(0^2,-12,-13,-14,-15-23,-25+34)$ & $\langle e_6,e_7 \rangle$\\
\specialrule{1pt}{0pt}{0pt}
\text{$23457$E} & $(0^2,-12,-13,-14-23,-15-24,-23)$ & $\langle e_6,e_7 \rangle$\\
\specialrule{1pt}{0pt}{0pt}
\text{$23457$F} & $(0^2,-12,-13,-14-23,-25+34,-23)$ & $\langle e_6, e_7 \rangle$\\
\specialrule{1pt}{0pt}{0pt}
\text{$23457$G} & $(0^2,-12,-13,-14-23,-15-24,-25+34)$ & $\langle e_6,e_7 \rangle$\\
\bottomrule[1pt]
\end{tabu}}
\end{table}
}

{\footnotesize 
\begin{table}[H]
\centering
\caption{7-dimensional indecomposable 6-step nilpotent Lie algebras}\label{7d-ind-6step-1}
\vspace{0.25 cm}
{\tabulinesep=1.2mm
\begin{tabu}{l|c|c}
\toprule[1.5pt]
NLA & structure equations & center\\
\specialrule{1pt}{0pt}{0pt}
\text{$123457$A} & $(0^2,-12,-13,-14,-15,-16)$ & $\langle  e_7 \rangle$\\
\specialrule{1pt}{0pt}{0pt}
\text{$123457$B} & $(0^2,-12,-13,-14,-15,-16-23)$ & $\langle  e_7 \rangle$\\
\specialrule{1pt}{0pt}{0pt}
\text{$123457$C} & $(0^2,-12,-13,-14,-15,-16-25+34)$ & $\langle  e_7 \rangle$\\
\specialrule{1pt}{0pt}{0pt}
\text{$123457$D} & $(0^2,-12,-13,-14,-15-23,-16-24)$ & $\langle  e_7 \rangle$\\
\specialrule{1pt}{0pt}{0pt}
\text{$123457$E} & $(0^2,-12,-13,-14,-15-23,-16-23-24)$ & $\langle  e_7 \rangle$\\
\specialrule{1pt}{0pt}{0pt}
\text{$123457$F} & $(0^2,-12,-13,-14,-15-23,-16-24-25+34)$ & $\langle  e_7 \rangle$\\
\specialrule{1pt}{0pt}{0pt}
\text{$123457$H} & $(0^2,-12,-13,-14-23,-15-24,-16-23-25)$ & $\langle  e_7 \rangle$\\
\specialrule{1pt}{0pt}{0pt}
\text{$123457\text{H}_1$} & $(0^2,-12,-13,-14-23,-15-24,16-23+25)$ & $\langle  e_7 \rangle$\\
\specialrule{1pt}{0pt}{0pt}
\text{$123457\text{I}(\lambda)$} & $(0^2,-12,-13,-14-23,-15-24,-16-\lambda \cdot 25+(\lambda-1)\cdot 34)$ & $\langle  e_7 \rangle$\\
\bottomrule[1pt]
\end{tabu}}
\end{table}
}

\noindent {\bf Conflict of Interest.} 
On behalf of both authors, the corresponding author states that there is no conflict of interest.

\noindent {\bf Data Availability Statement}
No additional data were generated or analyzed in this research.

\end{document}